\documentclass{sunil}
\usepackage[latin1]{inputenc}
\usepackage{amsmath}
\usepackage{amssymb}
\usepackage{latexsym}
\usepackage{verbatim}
\usepackage{graphics}
\usepackage{multirow}
\usepackage{amsfonts}
\usepackage{graphicx}
\usepackage{epstopdf}
\usepackage{pdflscape}
\usepackage{multicol}
\usepackage[margin=0.5in]{geometry}

\def\e{\varepsilon}

\begin{document}

\title{Invariants of Spatial Graphs}

\author{Blake Mellor}


\date{\today}

\maketitle

\section{Introduction}

The fundamental problem in knot theory is the classification of knots; in other words, the classification up to isotopy of all embeddings of $S^1$ in $S^3$.  Similarly, the fundamental problem in the theory of spatial graphs is, for any given graph $G$, the classification up to isotopy of the embeddings of $G$ in $S^3$.  In this sense, knot theory is simply a special case of the theory of spatial graphs, since a knot (resp. link) can be viewed as a cycle graph (resp. disjoint union of cycle graphs). 

\begin{figure} [ht]
    \begin{center}
    \scalebox{.7}{\includegraphics{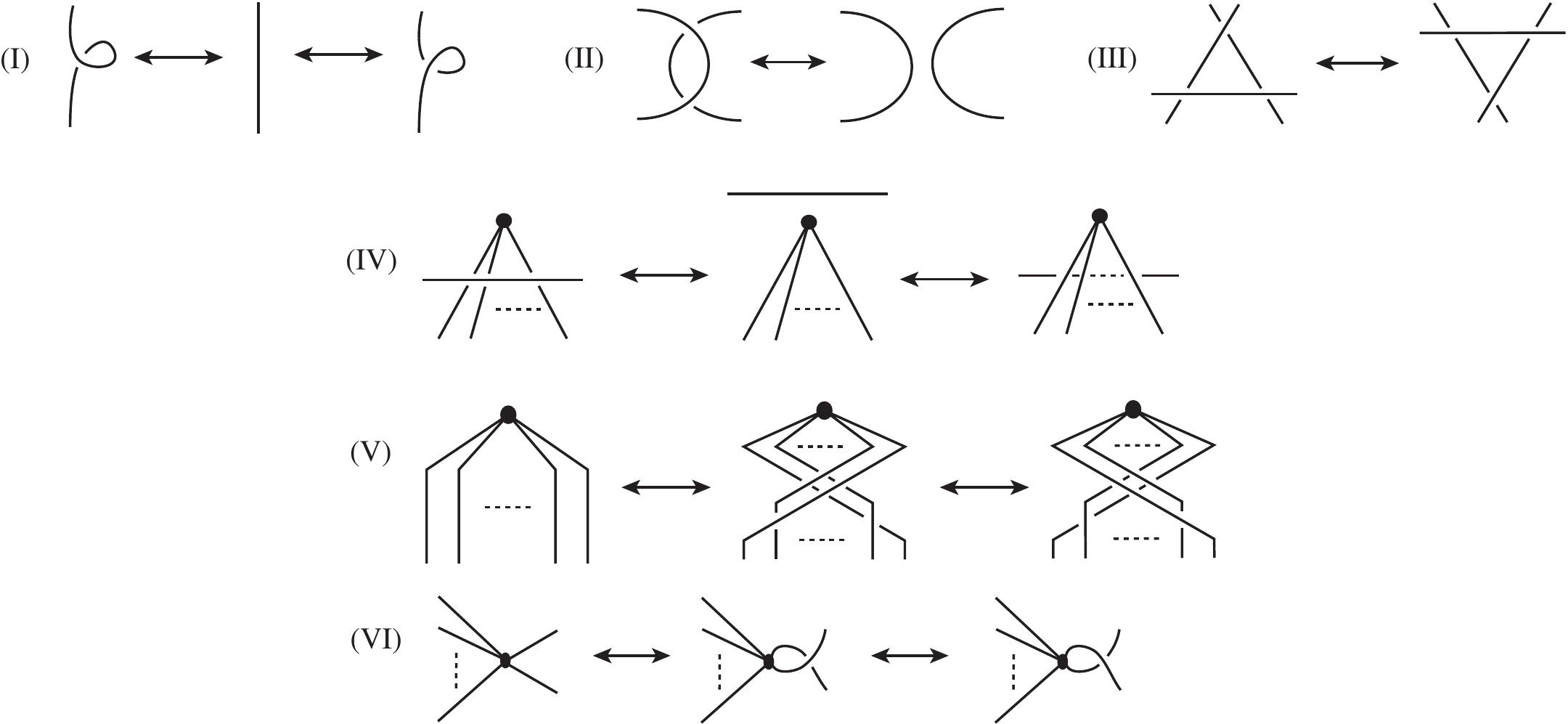}}
    \end{center}
    \caption{Reidemeister moves for spatial graphs} \label{F:reidemeister}
    \end{figure}

The two tasks in such a classification are determining when two embeddings are isotopic, and when they are not.  Both tasks are generally difficult. As with knots, two spatial graphs are isotopic if and only if their diagrams are equivalent modulo a set of {\em Reidemeister moves}, shown in Figure~\ref{F:reidemeister}; however, finding the specific sequence of moves between two equivalent diagrams can be extremely challenging. In fact, there are two slightly different notions of isotopic that are commonly considered.  Moves (I)-(V) in Figure \ref{F:reidemeister} generate {\em rigid (or flat) vertex isotopy}, in which the cyclic order of the edges at each vertex is fixed.  Including move (VI) generates {\em pliable vertex isotopy}, or simply {\em isotopy}, in which the order of the edges around a vertex may be changed.

Conversely, to determine that two embeddings are {\em not} isotopic requires an {\em invariant} - a function of the embeddings whose output is not changed by isotopies, and which takes different values on the two embeddings. Over the last century, many such invariants have been defined and studied for knots; in this article, we will look at some ways in which those ideas have been extended to other spatial graphs, with a particular emphasis on combinatorial and polynomial invariants.  In this article, we will mostly discuss invariants of pliable vertex isotopy, though we will discuss one important invariant of rigid vertex isotopy, the Yamada polynomial.

\section{Knots in graphs}

If a graph has no cycles (i.e. the graph is a tree), then all of its embeddings in $S^3$ are isotopic, so in any case of interest the graph will have cycles.  In the embedded graph, each of these cycles becomes a knot (possibly a trivial knot), and any collection of pairwise disjoint cycles becomes a link (possibly a trivial link). So to each embedding $G$ of a graph we can associate a collection $T(G)$ of knots and links (with multiplicity, as a particular knot type may appear many times).  More precisely, at each vertex choose two edges to connect, and delete the others; making this choice at every vertex leaves some number (possibly 0) of closed loops that form a link.  The set of all such links, for all choices of two edges at each vertex, is $T(G)$.  Kauffman \cite{kau} first observed that $T(G)$ is an isotopy invariant of $G$, and that many other (more easily comparable) invariants can be defined by applying our favorite knot and link invariants to the elements of $T(G)$. Conway and Gordon's famous result that every embedding of the complete graph on seven vertices contains a knot uses this perspective - they consider the sum of the Arf invariants of all Hamiltonian cycles in the embedded graph (mod 2) and show that this invariant is always nonzero \cite{cg}.

\begin{figure} [ht]
    $$G = \scalebox{.7}{\raisebox{-.4\height}{\includegraphics{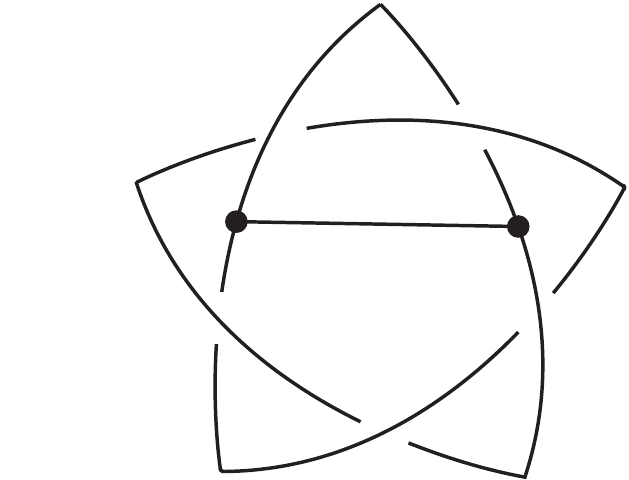}}} \qquad\longrightarrow \qquad T(G) = \left\{\scalebox{.7}{\raisebox{-.4\height}{\includegraphics{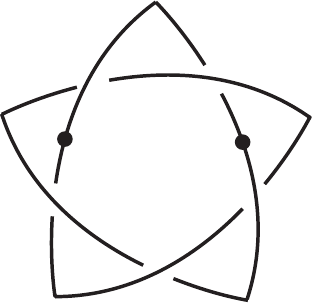}}}, \scalebox{.7}{\raisebox{-.4\height}{\includegraphics{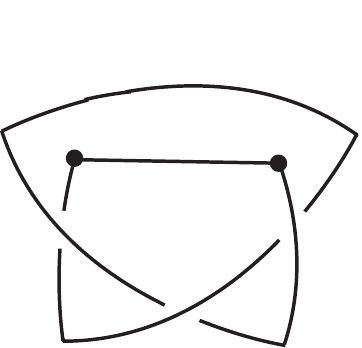}}}, \scalebox{.7}{\raisebox{-.4\height}{\includegraphics{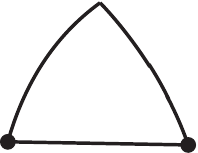}}} \right\}$$
    \caption{$T(G)$ for a $\theta$-graph} \label{F:T(G)}
    \end{figure}

As an example, consider the $\theta$-graph in Figure~\ref{F:T(G)} (this is $\theta$-graph $5_4$ in Moriuchi's table \cite{mo1}).  In this case $T(G)$ is the collection of {\em constituent knots} formed by selecting two of the three edges of the graph; these are a (2,5) torus knot, a trefoil knot, and an unknot. 

\section{The fundamental group and the Alexander polynomial of a spatial graph}

Perhaps the most important knot invariant is the {\em fundamental group} of the knot (i.e. the fundamental group of the knot complement).  The fundamental group can be used to define many other invariants, most famously the {\em Alexander polynomial} of the knot \cite{al}.  Similarly, the {\em fundamental group of a spatial graph $G$} is the fundamental group of the complement of the embedded graph, which is invariant under pliable vertex isotopy, and it can also be used to define other invariants.  An Alexander polynomial for spatial graphs was first defined by Kinoshita \cite{ki}. As with the Alexander polynomial of knots, the polynomial for spatial graphs arises from the first homology of the infinite cyclic cover of the graph complement.  However, also as with knots, it can be easily computed from the coefficient matrix for a system of linear relations that we can read off from the diagram of the spatial graph.

\begin{figure}[th]
\begin{center}
\scalebox{.8}{\includegraphics{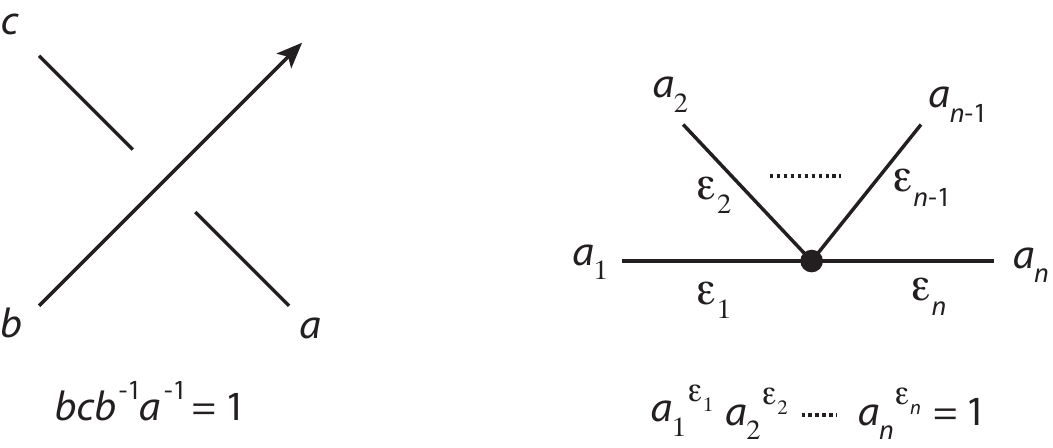}}
\end{center}
\caption{Wirtinger relations for $\pi_1(S^3 - G)$.}
\label{F:wirtinger}
\end{figure}

As with knots, the computation of the fundamental group for a spatial graphs depends on an {\em orientation} of the graph - i.e. an orientation of each edge of the graph.  For an oriented spatial graph $G$, the fundamental group $\pi_1(S^3-G)$ has a Wirtinger presentation constructed from a diagram for the graph, where the generators correspond to the arcs in the diagram.  The presentation has a relation at each crossing and vertex in the diagram, as shown in Figure~\ref{F:wirtinger}.  At a vertex, the {\it local sign} $\e_i$ of arc $a_i$ is $+1$ if the arc is directed into the vertex, and $-1$ if the arc is directed out from the vertex.

The Alexander polynomial is defined for a {\em balanced} oriented spatial graph.  This means each edge is given an integral {\em weight} so that the oriented sum of weights at each vertex is zero.  Then we get a linear relation at each crossing and vertex of a diagram of the graph, where the variables are the arcs of the diagram as shown in Figure \ref{F:wirtinger} (note that the orientation of the undercrossing arc on the left is not specified, since either orientation leads to the same relation).  In the crossing relation we will assume arc $b$ has weight $w_1$, and arcs $a$ and $c$ (on the same edge) have weight $w_2$.  In the vertex relation, arc $a_i$ has weight $w_i$.  Let $m_i = \e_1w_1 + \e_2w_2 + \cdots + \e_{i-1}w_{i-1} + \min\{\e_i, 0\}w_i$. Then the {\em Alexander relations} are:
\begin{itemize}
	\item[] Crossing relation: $(1-t^{w_2})b + t^{w_1}c - a = 0$ and 
	\item[] Vertex relation: $\sum_{i=1}^n{\e_i t^{m_i}a_i}  = 0$.
\end{itemize}
The Alexander matrix is the coefficient matrix for this system of linear relations.  For knots, there are only crossing relations, and the Alexander matrix is a square $r\times r$ matrix.  The Alexander polynomial is then any $(r-1) \times (r-1)$ minor of the matrix.  For graphs, however, the matrix is generally {\em not} square; rather, it is an $r \times s$ matrix with $s \geq r$.  In this case the Alexander polynomial is the greatest common divisor of the $(r-1) \times (r-1)$ minors of the matrix.  As for knots, the Alexander polynomial is well-defined only up to multiplication by $t^k$.

\begin{figure}[th]
\begin{center}
\scalebox{.5}{\includegraphics{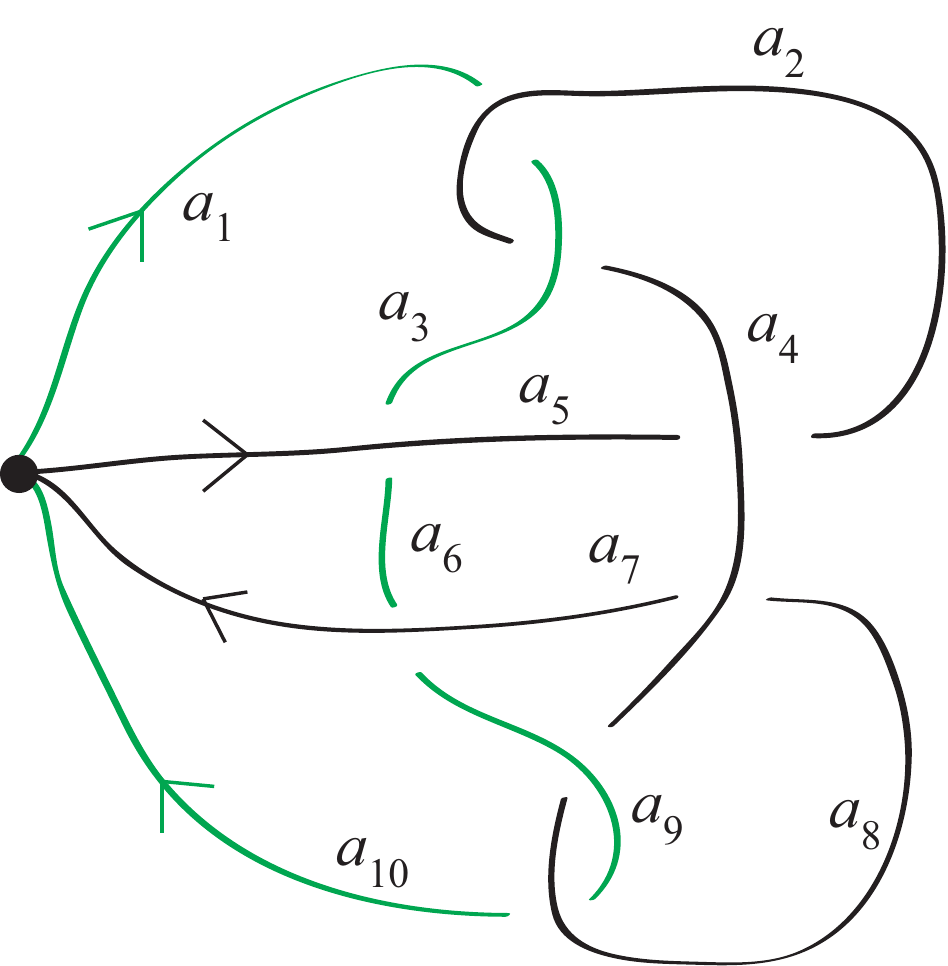}}
\end{center}
\caption{A spatial graph diagram $D$.}
\label{F:bouquet}
\end{figure}

As an example, consider the spatial graph $G$ with diagram $D$ shown in Figure \ref{F:bouquet}, with the edges given weights $x$ and $y$ (the edge containing arc $a_1$ has weight $x$).  Then the Alexander matrix is shown below (where the columns correspond to the arcs $a_1,\dots, a_{10}$ and the last row is the vertex relation):
$$\left[ \begin{matrix} -1 & 1-t^x & t^y & 0 & 0 & 0 & 0 & 0 & 0 & 0 \\ 0 & -1 & 1-t^y & t^x & 0 & 0 & 0 & 0 & 0 & 0\\ 0 & 0 & t^y & 0 & 1-t^x & -1 & 0 & 0 & 0 & 0 \\ 0 & t^y & 0 & 1-t^y & -1 & 0 & 0 & 0 & 0 & 0\\ 0 & 0 & 0 & 0 & 0 & -1 & 1-t^x & 0 & t^y & 0 \\ 0 & 0 & 0 & 1-t^y & 0 & 0 & -1 & t^y & 0 & 0 \\ 0 & 0 & 0 & t^x & 0 & 0 & 0 & -1 & 1-t^y & 0\\ 0 & 0 & 0 & 0 & 0 & 0 & 0 & 1-t^x & t^y & -1 \\ -t^{-x} & 0 & 0 & 0 & -t^{-x-y} & 0 & t^{-x-y} & 0 & 0 & t^{-x} \end{matrix}\right]$$
If both edges are given weight 1, then the Alexander polynomial is $\Delta(G)(t) = t^2 - 2t + 2$ (normalized so that the lowest term is the constant term).

\begin{figure}[htbp]
\begin{center}
\scalebox{.6}{\includegraphics{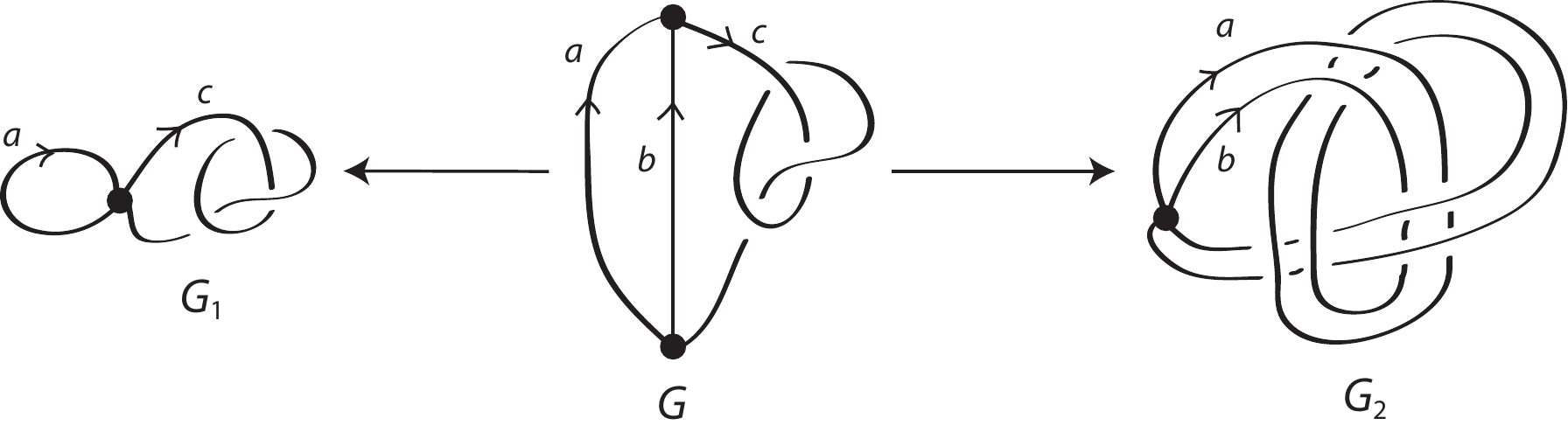}}
\end{center}
\caption{Contracting edges of a $\theta$-graph.}
\label{F:theta}
\end{figure}

Since the fundamental group of a spatial graph is a topological invariant of the exterior, we can contract a spanning tree of the graph without changing its fundamental group.  For example, the three graphs in Figure~\ref{F:theta} all have isomorphic fundamental groups, since $G_1$ and $G_2$ are each the result of contracting one edge of the $\theta$-graph $G$.  The Alexander polynomial, however, also depends on the weights on the edges.  For example, if the two arcs in $G_1$ are given weight 1, the Alexander polynomial is $\Delta(G_1)(t) = t^2-t+1$.  But if we give both arcs of $G_2$ weight 1, the Alexander polynomial is $\Delta(G_2)(t) = t^4 - t^2 +1$.  The reason for the difference is that these weightings result from different weightings on $G$; if we try to give all three arcs of $G$ a weight of 1, then the weighting is not balanced.  On the other hand, if we give arcs $a$ and $b$ weight 1 and arc $c$ weight 2, then $G$, $G_1$ and $G_2$ will all have an Alexander polynomial of $t^4-t^2+1$.

Further invariants can be derived from the Alexander matrix, such as the {\em determinant} of the spatial graph (which is calculated by substituting $-1$ for $t$ in the Alexander matrix, and then computing the greatest common divisor of the $(r-1) \times (r-1)$ minors), or whether the graph is {\em p-colorable} \cite{msw, mklp}.  Litherland \cite{lit} defined a more refined, and more powerful, Alexander polynomial for $\theta_n$-graphs (graphs with two vertices, and $n$ edges connecting them).  In his treatment, the edges of the graph were all oriented the same way, and there were no weights on the edges (though the edges are ordered).

While the Alexander polynomial of a spatial graph is constructed similarly to the Alexander polynomial of a knot, it is not as well understood.  For example, it is well-known which polynomials can be realized as the Alexander polynomial of a knot; for spatial graphs this question is still open.  Also, the Alexander polynomial satisfies a nice {\em skein relation} (particularly when normalized to give the {\em Alexander-Conway polynomial}) which provides an alternative method of computation.  We do not know of a similar skein relation for the Alexander polynomial of a spatial graph.

There are other algebraic invariants for spatial graphs as well. For example, Ghuman \cite{gh} defined a {\em longitude} of a cycle in a graph, and used them to define invariants analogous to Milnor's $\bar{\mu}$-invariants for links.  However, these are only isotopy invariants if the longitude is unique, which occurs when all vertices on the cycle have valence 3; where they are defined, however, they can contain more information than the fundamental group. Fleming \cite{fl} generalized Milnor's invariants in another way, looking at ``links" of spatial graphs with multiple connected components.

\section{The fundamental quandle of a spatial graph}

Another important algebraic object we can associate to a knot or spatial graph is the {\em fundamental quandle}.  The {\em fundamental quandle} of a knot was introduced in the early 1980's by Joyce \cite{joy}. Like the fundamental group, the fundamental quandle is defined by Wirtinger-type relations at the crossings.  However, Joyce proved that the fundamental quandle, unlike the fundamental group, is a complete knot invariant (of course, since nothing is free, quandles are much harder to compare than groups!). Niebrzydowski \cite{ni} extended the fundamental quandle to spatial graphs.  Though it seems unlikely that the fundamental quandle of a spatial graph is a complete invariant, as it is for knots, it can be used to define many useful invariants of a graph \cite{ii, jo, ni, os}.

A {\it  quandle} is a set $Q$ equipped with two binary operations $\rhd$ and $\rhd^{-1}$ that satisfy the following three axioms:
\begin{enumerate}
\item $x \rhd x =x$ for all $x \in Q$
\item $(x \rhd y) \rhd^{-1} y = x = (x \rhd^{-1} y) \rhd y$ for all $x, y \in Q$
\item $(x \rhd y) \rhd z = (x \rhd z) \rhd (y \rhd z)$ for all $x,y,z \in Q$
\end{enumerate}

The operation $\rhd$ is, in general, not associative. It is useful to adopt the exponential notation introduced by Fenn and Rourke in \cite{fr} and denote $x \rhd y$ as $x^y$ and $x \rhd^{-1} y$ as $x^{\bar y}$. With this notation, $x^{yz}$ will be taken to mean $(x^y)^z=(x \rhd y)\rhd z$ whereas $x^{y^z}$ will mean $x\rhd (y \rhd z)$. If $n$ is an integer, we will also let $x^{y^n} = x^{yy\cdots y}$, where the $y$ is repeated $n$ times.

Given a diagram for a spatial graph $G$, the fundamental quandle $Q(G)$ is the quandle whose generators are the arcs of the diagram, with additional Wirtinger-type relations at each crossing and vertex.  Using the labels on the arcs in Figure \ref{F:wirtinger}, the additional relations are:
\begin{enumerate}\setcounter{enumi}{3}
	\item $a^b = c$, and 
	\item For {\em every} generator $x$ of the quandle, $x^{a_1^{\e_1} \cdots a_n^{\e_n}} = x$, where $a_i^{\e_i} = a_i$ if $\e_i = +1$ and $a_i^{\e_i} = \overline{a_i}$ if $\e_i = -1$.
\end{enumerate}
Relation (4) is the same as for the fundamental quandle of a knot; relation (5) ensures invariance under Reidemeister moves (IV)-(VI), so the quandle is an invariant of pliable vertex isotopy.

Given a quandle $Q$, there is an {\em associated group} $As(Q)$, obtained by interpreting the quandle operation as conjugation (i.e. $x^y = y^{-1}xy$).  For a knot, the associated group of the fundamental quandle is isomorphic to the fundamental group.  However, for spatial graphs, these groups are generally not isomorphic (though there is always an epimorphism from the associated group of the fundamental quandle to the fundamental group).  In particular, while the abelianization of the fundamental group of a knot is always isomorphic to $\mathbb{Z}$, the abelianization of $As(Q(G))$ is equal to $\mathbb{Z}^E$, where $E$ is the number of edges in the graph.  So two cycles with different numbers of edges will be distinguished, even if their fundamental groups are isomorphic.

A particularly rich source of invariants derived from the fundamental quandle are {\em quandle colorings}.  A coloring of a diagram $D$ of a spatial graph $G$ by a quandle $X$ is an assignment of elements of $X$ to the arcs of $D$ so that relations (4) and (5) are satisfied at each crossing and vertex. Essentially, this is a quandle homomorphism from $Q(G)$ to $X$.  For a fixed $X$, the number of quandle colorings is an invariant \cite{ni}, and other coloring invariants of spatial graphs (such as those in \cite{iy, msw, mklp}) can be interpreted in terms of quandles.  These also lead to quandle homology and cohomology invariants \cite{ii}.  There is still much room to extend the work that has been done on knot quandles to spatial graphs.

\section{The Yamada polynomial}

In the 1980's, knot theorists discovered a family of new knot invariants, including the Jones, Kauffman and HOMFLY polynomials. Like the Alexander polynomial, these invariants satisfied nice skein relations which make them relatively easy to compute - several of them, however, are defined entirely in terms of this skein relation, rather than derived from a deeper topological invariant such as the fundamental group. The proof of invariance then relies on using the skein relation to show the value of the invariant is unchanged by Reidemeister moves.

This inspired attempts to construct similar, combinatorially defined, polynomial invariants for spatial graphs.  The best-known of these is the Yamada polynomial \cite{ya}, which can be defined as the unique polynomial $R(G)(A)$ which satisfies the following axioms:

\begin{enumerate}
    \item $R\left(\raisebox{-.3\height}{\includegraphics{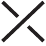}}\right) = AR\left(\raisebox{-.3\height}{\includegraphics{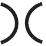}}\right)+A^{-1}R\left(\raisebox{-.3\height}{\includegraphics{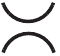}}\right) + R\left(\raisebox{-.3\height}{\includegraphics{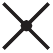}}\right)$
    \item $R(G) = R(G - e) + R(G/e)$, where $e$ is a nonloop edge in $G$, $G - e$ is the result of deleting $e$, and $G/e$ is the result of contracting $e$.
    \item $R(G_1 \amalg G_2) = R(G_1)R(G_2)$, where $\amalg$ denotes disjoint union.
    \item $R(G_1 \vee G_2) = -R(G_1)R(G_2)$, where $G_1 \vee G_2$ is the graph obtained by joining $G_1$ and $G_2$ at any single vertex.
    \item $R(B_n) = -(-\sigma)^n$, where $B_n$ is the $n$-leafed bouquet of circles and $\sigma = A + 1 + A^{-1}$.  In particular, if $G$ is a single vertex, $R(G) = R(B_0) = -1$.
    \item $R(\emptyset) = 1$
\end{enumerate}

Using these skein relations, $R(G)$ can be computed by reducing $G$ to combinations of bouquets of circles.  $R(G)$ is an invariant of spatial graphs up to {\it regular rigid vertex isotopy}, meaning that it is invariant under moves (II), (III) and (IV) in Figure \ref{F:reidemeister}, but not moves (I), (V) or (VI).  The behavior of $R$ under these moves is given by the following formulas:

\begin{enumerate}
	\item $R\left(\scalebox{.3}{\raisebox{-.3\height}{\includegraphics{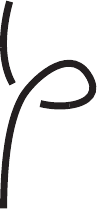}}}\right) = A^2R\left(\vert\right)$ \quad and \quad $R\left(\scalebox{.4}{\raisebox{-.3\height}{\includegraphics{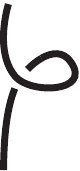}}}\right) = A^{-2}R\left(\vert\right)$
	\item $R\left(\scalebox{.3}{\raisebox{-.3\height}{\includegraphics{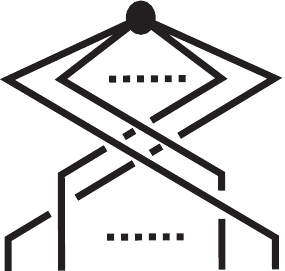}}}\right) = (-A)^nR\left(\scalebox{.3}{\raisebox{-.3\height}{\includegraphics{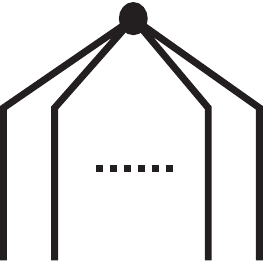}}}\right)$ \quad and \quad $R\left(\scalebox{.3}{\raisebox{-.3\height}{\includegraphics{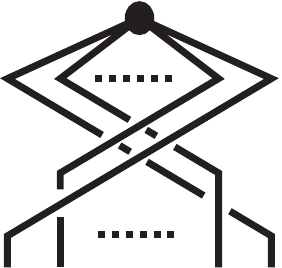}}}\right) = (-A)^{-n}R\left(\scalebox{.3}{\raisebox{-.3\height}{\includegraphics{R5a-eps-converted-to.pdf}}}\right)$
	\item $R\left(\scalebox{.3}{\raisebox{-.3\height}{\includegraphics{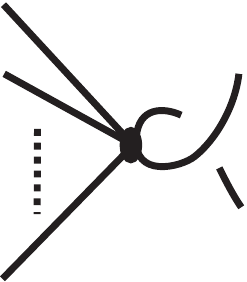}}}\right) = -AR\left(\scalebox{.3}{\raisebox{-.3\height}{\includegraphics{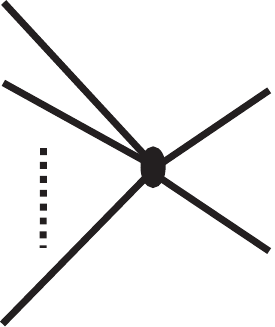}}}\right) - (A^2+A)R\left(\scalebox{.3}{\raisebox{-.3\height}{\includegraphics{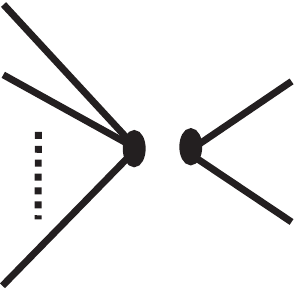}}}\right)$  
	\item $R\left(\scalebox{.3}{\raisebox{-.3\height}{\includegraphics{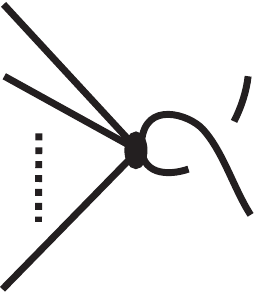}}}\right) = -A^{-1}R\left(\scalebox{.3}{\raisebox{-.3\height}{\includegraphics{R6d-eps-converted-to.pdf}}}\right) - (A^{-2}+A^{-1})R\left(\scalebox{.3}{\raisebox{-.3\height}{\includegraphics{R6c-eps-converted-to.pdf}}}\right)$
\end{enumerate}

We can obtain an invariant of rigid vertex isotopy (invariance under moves (I) - (V)) by defining $\overline{R}(G) = (-A)^{-m}R(G)$, where $m$ is the smallest power of $A$ in $R(G)$.  This will still not be invariant under move (VI), however, so it is not an invariant of pliable vertex isotopy.  The exception is when the maximum degree of the vertices of the graph is 3 or less, since a move of type (VI) on a vertex of degree 3 is equivalent to a move of type (V) followed by moves of types (IV) and (I).  So if the maximum degree of the vertices is 3 or less, $\overline{R}(G)$ is an invariant of pliable vertex isotopy. If the graph $G$ is simply a link, then $R(G)$ is a specialization of the Dubrovnik version of the Kauffman polynomial; if $G$ is a knot, then $R(G)$ is the Jones polynomial of the $(2,0)$-cabling of the knot \cite{ya}.

\begin{figure} [ht]
    $$\scalebox{.7}{\raisebox{-.4\height}{\includegraphics{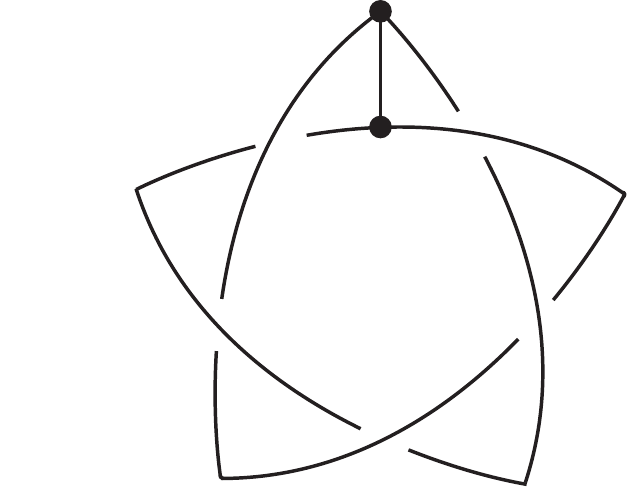}}}$$
    $$\overline{R}(G) = -1-A-A^2-A^3-A^4-A^{10}-A^{12}-A^{14}+A^{16}+A^{18}$$
    $$\scalebox{.7}{\raisebox{-.4\height}{\includegraphics{theta_5_4-eps-converted-to.pdf}}}$$
    $$\overline{R}(G) = -1-A-A^2-A^3-2A^4-A^5-A^6-A^7+A^9+A^{11}+A^{13}+A^{16}-A^{17}$$
    \caption{Yamada polynomials for $\theta$-graphs $5_3$ and $5_4$} \label{F:yamada}
    \end{figure}

The skein relations make the Yamada polynomial convenient to work with and relatively easy to compute (with the help of a computer); it is probably the most popular invariant of spatial graphs.  It has been used to provide necessary conditions for a spatial graph to be {\em p-periodic} (i.e. symmetric under an action of $\mathbb{Z}_p$ on $S^3$ by homeomorphisms) \cite{ch, ma}.  The Yamada polynomial has also been used to provide a lower bound on the number of crossings in a spatial graph \cite{mot}.  The Yamada polynomial is often used to help develop ``knot tables" for particular graphs, such as $\theta$-graphs, handcuff graphs and bouquet graphs \cite{dv, mo1, mo2, oy}.  For example, Figure~\ref{F:yamada} shows the $\theta$-graphs $5_3$ and $5_4$ from Moriuchi's table \cite{mo1}, along with their Yamada polynomials $\overline{R}(G)$.  (These graphs are also distinguished by their constituent knots, since $T(5_3)$ contains only one knot, while $T(5_4)$ contains two.) The Yamada polynomial is important enough that is has become an object of study in its own right; for example, Li, Lei, Li and Vesnin \cite{lllv} have studied the distribution of zeros of the Yamada polynomials of certain classes of spatial graphs.

The fact that the Yamada polynomial is, in general, only an invariant of rigid vertex isotopy is its primary drawback.  Other polynomials are invariants of pliable vertex isotopy, but these are generally harder to compute.  In addition to the Alexander polynomial discussed in the last section, the Yokota polynomial \cite{yo} is an isotopy invariant which agrees with the Yamada polynomial for graphs of constant valence 3; however, its computation is far more complex, involving transforming the graph into a large linear combination of link diagrams, and summing invariants of all the terms.  The Thompson polynomial \cite{th} is also an isotopy invariant, which can be computed recursively, but each stage of the recursion requires determining topological facts about the boundary of a three-manifold, such as whether it is compressible.

\section{Other invariants of spatial graphs} 

We have only scratched the surface of spatial graph invariants.  We have discussed some invariants that are defined for all (or at least most) graphs, but there are many others defined for smaller classes of graphs that are still useful.  For example, the Simon invariants were initially defined only for embeddings of the complete graphs $K_5$ and $K_{3,3}$, and yet they have turned out to be surprisingly useful \cite{ffn, nt, st}. Or, we can look at invariants for equivalence relations other than rigid vertex isotopy and pliable vertex isotopy. Many such relations have been studied for knots and links, and nearly all of these have an analogue (or even several analogues) for spatial graphs. For instance, Taniyama introduced homotopy and homology equivalence relations on spatial graphs \cite{ta1, ta2, ta3}, and many others have studied these ideas as well \cite{fn, st}. Other ideas from knot theory, such as {\em finite type invariants}, can also be extended to spatial graphs \cite{hj, ty}.  We could go on, but it is enough to point out that invariants of spatial graphs have been studied for only 30 years, while invariants of knots and links have been studied for over a century; there is much to be done!

\small

\normalsize


\begin{thebibliography}{10}

\bibitem{al} J.W. Alexander, Topological invariants of knots and links, {\it Trans. Amer. Math. Soc.} (1928), pp. 275-306

\bibitem{ch} N. Chbili, Skein algebras of the solid torus and symmetric spatial graphs, {\it Fundamenta Math.} v. 190 (2006)

\bibitem{cg} J. Conway and C. Gordon,  Knots and links in spatial graphs, {\it J. of Graph Theory} v. 7 (1983), pp. 445-453.

\bibitem{dv} A. Dobrynin and A. Vesnin, The Yamada polynomial for graphs, embedded knot-wise into three-dimensional space, {\em Vychisl. Sistemy} v.155 (1996), pp. 37-86 (in Russian). An English translation is available at https://www.researchgate.net/publication/266336562.

\bibitem{fr} R. Fenn and C. Rourke, Racks and links in codimension two, {\em J. of Knot Theory and Its Ramif.} v. 1 (1992), pp. 343-406

\bibitem{ffn} E. Flapan, W. Fletcher and R. Nikkuni, Reduced Wu and generalized Simon invariants for spatial graphs, {\em Math. Proc. Camb. Phil. Soc.} v. 156 (2014), pp. 521-544

\bibitem{fl} T. Fleming, Milnor invariants for spatial graphs, {\em Topology Appl.}, v. 155 (2008), pp. 1297-1305

\bibitem{fn} T. Fleming and R. Nikkuni, Homotopy on spatial graphs and the Sato-Levine invariant, {\em Trans. Amer. Math. Soc.} v. 361 (2009), pp. 1885-1902

\bibitem{gh} S. Ghuman, Invariants of graphs, {\it J. Knot Theory Ramif.} v. 9 (2000), pp. 31-92

\bibitem{hj} Y. Huh and G.T. Jin,  $\theta$-curve polynomials and finite-type invariants, {\em J. Knot Theory Ramif.} v. 11 (2002), pp. 555-564

\bibitem{ii} A. Ishii and M. Iwakiri, Quandle cocycle invariants for spatial graphs and knotted handlebodies, {\em Canad. J. Math.} v. 64 (2012), pp. 102-122

\bibitem{iy} Y. Ishii and A. Yasuhara, Color invariant for spatial graphs, {\it J. Knot Theory Ramif.}, v. 6 (1997), pp. 319-325

\bibitem{jo} Y. Jang and K. Oshiro, Symmetric quandle coloring for spatial graphs and handlebody-links, {\em J. Knot Theory Ramif.} v. 21 (2012)

\bibitem{joy} D. Joyce, A classifying invariant of knots, the knot quandle, {\em J. of Pure and Applied Algebra} v. 23 (1982), pp. 37-65

\bibitem{kau} L.H. Kauffman, Invariants of graphs in three-space, {\it Trans. Amer. Math. Soc.} (1989), pp. 697-710

\bibitem{ki} S. Kinoshita, Alexander Polynomials as Isotopy Invariants I, {\it Osaka Math. J.} v. 10 (1958), pp. 263-271

\bibitem{lllv} M. Li, F. Lei, F. Li and A. Vesnin, On Yamada polynomial of spatial graphs obtained by edge replacements, preprint (2018).  Available at arXiv:1801.09075

\bibitem{lit} R. Litherland, The Alexander module of a knotted theta-curve, {\it Math. Proc. Camb. Phil. Soc.} v. 106 (1989), pp. 95-106

\bibitem{ma} Y. Marui, The Yamada polynomial of spatial graphs with $\mathbb{Z}_n$-symmetry, {\it Kobe J. Math.} v. 18 (2001), pp. 23-49

\bibitem{msw} J. McAtee, D. Silver and S. Williams, Coloring spatial graphs, {\it J. Knot Theory Ramif.} v. 10 (2001), pp. 109-120

\bibitem{mklp} B. Mellor, T. Kong, A. Lewald and V. Pigrish, Colorings, determinants and Alexander polynomials for spatial graphs, {\it J. Knot Theory Ramif.} v. 25 (2016)

\bibitem{mo1} H. Moriuchi, An enumeration of theta-curves with up to seven crossings, {\it J. Knot Theory Ramif.} v. 18 (2009), pp. 167-197

\bibitem{mo2} H. Moriuchi, A table of $\theta$-curves and handcuff graphs with up to seven crossings, in {\it Noncommutativity and singularities}, Adv. Stud. Pure Math., v. 55, Math. Soc. Japan, Tokyo, 2009, pp. 281-290

\bibitem{mot} T. Motohashi, Y. Ohyama and K. Taniyama, Yamada polynomial and crossing number of spatial graphs, {\em Rev. Mat. Univ. Complut. Madrid} v. 7 (1994), pp. 247-277

\bibitem{ni} M. Niebrzydowski, Coloring invariants of spatial graphs, {\it J. Knot Theory Ramif.} v. 19 (2010), pp. 829-841

\bibitem{nt} R. Nikkuni and K. Taniyama, Symmetries of spatial graphs and Simon invariants, {\em Fund. Math.} v. 205 (2009), pp. 219-236

\bibitem{os} K. Oshiro, On pallets for Fox colorings of spatial graphs, {\it Topology and its Applications} v. 159 (2012), pp. 1092-1105

\bibitem{oy} N. Oyamaguchi, Enumeration of spatial 2-bouquet graphs up to flat vertex isotopy, {\it Topology and its Appl.} v. 196 (2015), pp. 805-814

\bibitem{st} R. Shinjo and K. Taniyama, Homology classification of spatial graphs by linking numbers and Simon invariants, {\em Topology Appl.} v. 134 (2003), pp. 53-67

\bibitem{ta1} K. Taniyama, Link homotopy invariants of graphs in $\mathbb{R}^3$, {\em Rev. Mat. Univ. Complut. Madrid} v. 7 (1994), pp. 129-144.

\bibitem{ta2} K. Taniyama, Cobordism, homotopy and homology of graphs in $\mathbb{R}^3$, {\em Topology} v. 33 (1994), pp. 509-523.

\bibitem{ta3} K. Taniyama, Homology classification of spatial embeddings of a graph, {\em Topology Appl.} v. 65 (1995), pp. 205-228

\bibitem{ty} K. Taniyama and A. Yasuhara, Local moves on spatial graphs and finite type invariants, {\em Pacific J. Math.} v. 211 (2003), pp. 183-200

\bibitem{th} A. Thompson, A polynomial invariant of graphs in 3-manifolds, {\it Topology} v. 31 (1992), pp. 657-665

\bibitem{ya} S. Yamada, An invariant of spatial graphs, {\it J. Graph Theory} v. 13 (1989), pp. 537-551

\bibitem{yo} Y. Yokota, Topological invariants of graphs in 3-space, {\it Topology} v. 35 (1996), pp. 77-87

\end{thebibliography}
\end{document}